\documentclass[11pt]{article}%\documentclass[10pt]{article}%
\usepackage{amsmath,latexsym,epsf}
\begin{document}

%% ----------------------------------------------------------------------
\newcommand{\nc}{\newcommand}
\newcommand{\mc}{\mathcal}
\def\PP#1{{\mathrm{Pres}}^{#1}{T}\setcounter{equation}{0}}
\def\ns{$n$-star}\setcounter{equation}{0}
\def\nt{$n$-tilting}\setcounter{equation}{0}
\def\Ht#1{{{\mathrm{Hom}}_R(T,{#1})}\setcounter{equation}{0}}
\def\qp#1{{${(#1)}$-quasi-projective}\setcounter{equation}{0}}
\def\mr#1{{{\mathrm{#1}}}\setcounter{equation}{0}}
%%%%%%%%FROM latexexam.tex
%\theoremstyle{definition}
\newtheorem{Th}{Theorem}[section]%[section]
\newtheorem{Def}[Th]{Definition}
\newtheorem{Lem}[Th]{Lemma}
\newtheorem{Pro}[Th]{Proposition}
\newtheorem{Cor}[Th]{Corollary}
\newtheorem{Rem}[Th]{Remark}
\newtheorem{Exm}[Th]{Example}
\def\Pf#1{{\noindent\bf Proof}.\setcounter{equation}{0}}
\def\>#1{{ $\Rightarrow$ }\setcounter{equation}{0}}

\def\bskip#1{{\vskip 15pt}}%\setcounter{equation}{0}}
\def\zskip#1{{\vskip 10pt}}%\setcounter{equation}{0}}
\def\sskip#1{{\vskip 5pt}}%\setcounter{equation}{0}}

\def\bg#1{\begin{#1}\setcounter{equation}{0}}
\def\ed#1{\end{#1}\setcounter{equation}{0}}
\def\TFAE{The following are equivalent:}
%%%%%%%%%%

%%%%%%%%%%%%%%%%%%%%%%%%%%%%%%%%%%%%%%%%%%%%%%%%%%%%%%%%%%%%%%%%%%%%%%%%%%%%%%%%%
%**************************标题、摘要、分类号、关键字**************************

\title{\bf Auslander Bounds and Homological Conjectures}
\smallskip
\author{\small {Jiaqun Wei{\thanks {Supported by the National Natural Science Foundation of China (Grant No. 10971099)}}}\\
%\small Department of Mathematics,
%%
%Nanjing Normal University \\
%%
%\small Nanjing 210097,
%%
%P.R.China\\ \small Email: weijiaqun@njnu.edu.cn
}
\date{}
\maketitle
\baselineskip 18pt%17pt%16pt%\baselineskip  25.5pt %
%
% Abstract ------------------------------------------------------
%
\begin{abstract}
\baselineskip 14pt%16pt%17pt%16pt%\baselineskip  14pt%25.5pt %

Inspired by recent works on rings satisfying Auslander's
conjecture, we study invariants, which we call Auslander bounds,
and prove that they have strong relations to some homological
conjectures.

\vskip 10pt


\end{abstract}
%\smallskip
%
\vskip 30pt
% ----------------------------------------------------------------------
% ----------------------------------------------------------------------
%\def\baselinestretch{1}
%
\baselineskip 16pt%18pt%17pt%\baselineskip  14pt%25.5pt %

\section{Introduction}
%{\noindent \Large\bf Introduction}
%
%\vskip 15pt
%
%
\hskip 15pt
%
%%% ----------------------------------------------------------------------

Throughout this paper, rings are associative with nonzero
identities and modules are left modules unless otherwise
specified. Let $R$ be a ring, we denote by $\mr{Mod}R$ the
category of all left $R$-modules, and we denote by $\mr{mod}R$ the
full subcategory of all $R$-modules having finitely generated
projective resolutions.

%\sskip\

Auslander posed the following conjecture (cf. [\ref{H}]).

\bg{verse}

(Auslander Conjecture): {\it Let $R$ be an artin algebra. Then for
every $M\in \mr{mod}R$, there exists an integer $b_M$ such that
$\mr{Ext}_R^i(M,N)=0$ for all $i>b_M$ and for all $R$-modules
$N\in \mr{mod}R$ satisfying $\mr{Ext}_R^j(M,N)=0$ for all
sufficiently large $j$.}

\ed{verse}

It is known that if the Auslander conjecture holds for all finite
dimensional algebras then the finitistic dimension conjecture is
true for all finite dimensional algebras [\ref{H}]. However, the
Auslander conjecture fails in general by counterexamples in
[\ref{JS}, \ref{S}]. Rings satisfying the assertion in the
Auslander conjecture are studied in [\ref{CH}, \ref{HJ}, \ref{M}].
In [\ref{CH}] the authors investigate in detail the relationship
between such rings and some homological conjectures, for instance,
the Auslander-Reiten conjecture, which we recall as follows.

\bg{verse} Auslander-Reiten Conjecture: {\it Let $R$ be a ring. If
$M\in\mr{mod}R$ and $\mr{Ext}_R^i(M,M\oplus R)=0$ for all $i>0$,
then $M$ is projective.} \ed{verse}

%In particular, the following results are
%proved.
%
%\vskip 10pt
%
%\noindent{\bf [\ref{CH}, Theorem A]} Let $R$ be a left-noetherian
%ring satisfying (AC), and let $M\in \mr{mod}R$. If
%$\mr{Ext}_R^i(M,M) = 0$ for all sufficiently large $i$ and
%$\mr{Ext}_R^i(M,R) = 0$ for all $i>0$, then $M$ is projective.
%
%\vskip 10pt
%
%\noindent{\bf [\ref{CH}, Theorem B]} Let $R$ be a two-sided
%noetherian ring. If $R$ and $R^o$ satisfy (AC) and (1) $R$ is an
%artin algebra, or (2) $R$ has a dualizing complex, then
%$\mr{id}R<\infty$ if and only if $\mr{id}_{R^o}R<\infty$ (whence,
%$\mr{id}R = \mr{id}_{R^o}R$ by [\ref{I}]).
%
%
%\vskip 10pt

%and the relations between Auslander bounds and some homological
%conjectures

In this paper, we continue the study and focus on the number
$b_M$, which we called the left Auslander bound of $M$ and denote
by $\mr{lAb}M$, for the fixed module $M$. Note that the Auslander
bound can be $\infty$. We prove the following result.

%\sskip\

%\vskip 5pt

\bg{verse}{\bf Theorem} If $M\in\mr{mod}R$ and
$\mr{Ext}_R^i(M,M\oplus R)=0$ for all sufficiently large $i$, then
$\mr{lAb}M$ coincides with the projective dimension of
$M$.\ed{verse}

%\sskip\

The result has clear relation with the Auslander-Reiten
conjecture. More generally, it relates Auslander bounds with an
equivalent version of the Wakamatsu-tilting conjecture (EWTC for
short, see Section 3 for details), which asserts that an
$R$-module $T\in\mr{mod}R$ is tilting if and only if (1)
$\mr{Ext}_R^i(T,T)=0$ for all $i>0$ and (2) there is an exact
sequence $0\to R\to T_0\to\cdots\to T_n\to 0$ for some $n$, where
each $T_i\in\mr{add}T$. We refer to [\ref{AHK}] for the history
and development of tilting theory. In fact, we obtain that if $R$
satisfies the condition that $\mr{lAb}M<\infty$ for every
$M\in\mr{mod}R$, then $R$ satisfies the conjecture (EWTC). This
extends [\ref{CH}, Theorem A]. As we see, all rings with finite
global repetition index, in particular, the finite dimensional
algebra $\mc{O}/\pi\mc{O}$, where $\mc{O}$ is a classical order of
finite global dimension over a discrete valuation ring $D$ with
uniformizing parameter $\pi$ and residue class field $K$,
satisfies the condition in the last result.

The above theorem is proved in Section 2, after some
investigations on basic properties of Auslander bounds. Relations
between homological conjectures and Auslander bounds are presented
in Section 3, where we also formulate some new homological
conjectures.

%We investigate
%the basic properties of Auslander bounds and give some methods to
%calculate them (Proposition \ref{GD} and \ref{syg}, Theorems
%\ref{mth} and \ref{proself}). As we show in the paper, the results
%obtained in this way are also useful to give a new sight of some
%known theorems and to obtain new results (Corollary \ref{Cor-id},
%Propositions \ref{Ewtch}, \ref{gsc} and \ref{uc}).
%%
% %or moreover, to improve them (e.g., [\ref{CH},
%%Theorems A and B]) by removing the non-necessary conditions
%%
%
%Since the Auslander Conjecture fails in general even for artin
%algebras, we naturally consider the (left version of) finitistic
%Auslander Conjecture as follows.
%
%
%\vskip 10pt
%
%(fAC) {\it For every artin algebra $R$, there is an integer $n$
%such that the left Auslander bound of every finitely generated
%$R$-module with finite left Auslander bound is not more than $n$.
%}
%
%\vskip 10pt
%
%We show that the finitistic Auslander bound conjecture for all
%artin algebras implies the little finitistic dimension conjecture
%for all artin algebras (Proposition \ref{fac-fdc}) and give some
%cases in which the finitistic Auslander bound conjecture holds
%(Propositions \ref{id-fac}).

\bskip\

We introduce some notions in the following.

%Given an $R$-module $M$ and an integer $t\ge 0$, we denote by
%$M^{>t}$ ($M_f^{>t}$, resp.) the subcategory of all modules $N$
%($N\in \mr{mod}R$, resp.) such that $\mr{Ext}_R^i(M,N)=0$ for all
%$i>t$. The notion $M^{\gg 0}$ ($M_f^{\gg 0}$, resp.) denotes the
%subcategory of all modules $N$ such that $N\in M^{>t}$ ($N\in
%M_f^{>t}$, resp.) for some $t>0$.
%
Let $R$ be a ring and $C,D\in\mr{Mod}R$. Let $t$ be a non-negative
integer. By $\mr{Ext}_R^{>t}(C,D)=0$ we mean that
$\mr{Ext}_R^{i}(C,D)=0$ for all $i>t$. By
$\mr{Ext}_R^{\gg}(C,D)=0$ we mean that $\mr{Ext}_R^{>t}(C,D)=0$
for some $t$. Given an $R$-module $M$ and an integer $t\ge 0$, we
denote by $M^{>t}$ ($^{>t}M$, resp.) the subcategory of all
modules $N$ such that $\mr{Ext}_R^{>t}(M,N)=0$
($\mr{Ext}_R^{>t}(N,M)=0$, resp.). The notions $M^{\gg}$ and
$^{\gg}M$ are defined similarly.

% The notion $M^{\gg}$
%($^{\gg}M$, resp.) denotes the subcategory of all modules $N$ such
%that $N\in M^{>t}$ ($N\in {^{>t}M}$, resp.) for some $t\ge 0$.

% We use the notation
%$\mr{Ext}_R^{\gg 0}(\mc{M},\mc{N})$ to denote that  there is some
%$t$ such that  $\mr{Ext}_R^{>m}(\mc{M},\mc{N})$

For an $R$-module $M$, we denote by $\mr{add}M$ the class of all
modules isomorphic to direct summands of finite direct sums of
copies of $T$. We use $\mr{pd}M$ ($\mr{id}M$, $\mr{fd}M$, resp.)
to denote projective (injective, flat, resp.) dimension of $M$.

We denote by $R^o$ the opposite ring of $R$. Thus $\mr{Mod}R^o$ is
the category of all right $R$-modules. In case that $R$ is an
artin algebra, we denote by $\mathbf{D}$ the usual dual functor
between $\mr{mod}R$ and $\mr{mod}R^o$. %Thus, $(\mathbf{D}(M))_R$
%is just the dual module of $M\in\mr{mod}R$.

%
%%
%\vskip 15pt%
% \noindent{\bf Main Notions} Let $R$ be a ring. Let $\mc{C}$
%be a class of $R$-modules, we denote
%$\mc{C}_f=\mc{C}\bigcap\mr{mod}R$. If $\mr{X}M$ defines an
%invariance of the $R$-module $M$ in terms of $\mr{Mod}R$, then we
%denote by $\mr{X}_fM$ the invariance in the same way but in terms
%of $\mr{mod}R$.

%(i.e.,considered in the realm of entirely finitely generated modules)

%\vskip 15pt
%
%For example, $M_f^{>t}=M^{>t}\bigcap\mr{mod}R$ and
%

%%% ----------------------------------------------------------------------
\vskip 30pt
%---------------------------------------------------------------------------------------------

\section{Auslander bounds}
%%
%\hskip 15pt
%%
%%%% ----------------------------------------------------------------------

%\vskip 15pt\ \

\hskip 15pt%\sskip\

Throughout this section, we fix $R$ a ring. If a class $\mc{C}$
consists of a single $R$-module, say $C$, then we use $C$ instead
of $\mc{C}$.

We introduce the following notion.

\vskip 15pt

\bg{Def}\label{ABdef}%
Let $\mc{C},\mc{D}$ be two classes of $R$-modules. The Auslander
bound of the pair $(\mc{C},\mc{D})$, denoted by
$\mr{Ab}(\mc{C},\mc{D})$, is defined in the following way:

\begin{verse}
 If
there are no $C\in\mc{C}$ and $D\in\mc{D}$ such that
$\mr{Ext}^{\gg}_R(C,D)=0$, then $\mr{Ab}(\mc{C},\mc{D})=-1$;\\
Otherwise, $\mr{Ab}(\mc{C},\mc{D})$ is the minimal non-negative
integer $m$ such that $\mr{Ext}^{>m}_R(C,D)=0$ for any
$C\in\mc{C}$ and $D\in\mc{D}$ with $\mr{Ext}^{\gg}_R(C,D)=0$, or
$\infty$ if no such minimal integer exists.
%Let $M\in \mr{Mod}R$. The left $($right, resp.$)$
%$\mc{C}$-Auslander bound of $M$, denoted by $\mr{lAb}_{\mc{C}}M$
%$(\mr{rAb}_{\mc{C}}M$, resp.$)$,  is defined to be the minimal
%non-negative integer $m$ such that $M^{\gg 0}\bigcap
%\mc{C}=M^{>m}\bigcap \mc{C}$\ \ $( {^{\gg 0}M}\bigcap
%\mc{C}={^{>m}M}\bigcap \mc{C}$, resp.$)$, or $\infty$, if no such
%integer exists.
\end{verse}
\ed{Def}%

 It is easy to see
that $\mr{Ab}(C,D)<\infty$ for any two $R$-modules $C,D$.

Let $M\in\mr{Mod}R$. It is also easy to see that
$\mr{Ab}(M,{\mc{C}})$ is just the minimal non-negative integer $m$
such that $M^{\gg}\bigcap \mc{C}=M^{>m}\bigcap \mc{C}$, or
$\infty$ if no such integer exists, or $-1$ if $M^{\gg}\bigcap
\mc{C}=\O$. Similarly $\mr{Ab}({\mc{C}},M)$ is just the minimal
non-negative integer $m$ such that $\mc{C}\bigcap {^{\gg}M}
=\mc{C}\bigcap {^{>m}M}$, or $\infty$ if no such integer exists,
or $-1$ if $\mc{C}\bigcap {^{\gg}M}=\O$.

%%%%%%%%%%%%%%%%%%%%%%%%%%%%%%%%%%%%%%%%%%%%%%%%%%%%%%%%%%%%%%%%%%%%%%%%%%%%%%%%%%
%\vskip 15pt
%\sskip \ and $\mc{C}$ be a class of $R$-modules

 We use the following simple notions for $M\in\mr{Mod}R$.

\begin{verse}$\bullet\ \ \ \mr{LAb}M:=\mr{Ab}(M,\mr{Mod}R)$\ \ \ \ (called the big left Auslander bound of $M$);\\
$\bullet\ \ \ \mr{lAb}M\ :=\mr{Ab}(M,\mr{mod}R)$ \ \ \ \ (the small left Auslander bound of $M$);\\
$\bullet\ \ \ \mr{RAb}M:=\mr{Ab}(\mr{Mod}R,M)$ \ \ \ \ (the big right Auslander bound of $M$);\\
$\bullet\ \ \ \mr{rAb}M\ :=\mr{Ab}(\mr{mod}R,M)$ \ \ \ \ (the small right Auslander bound of $M$).%
%
%
%$\bullet\ \ \ \mr{rAb}M:=\mr{Ab}(\mr{mod}R,M)$

\end{verse}

It is easy to see that the above Auslander bounds are
non-negative.

%
%In case $\mc{C}=\mr{Mod}R$, we simply denote it by $\mr{lAb}M$
%(resp. $\mr{rAb}M$).

%We denote by $\mr{lAb}|_{\mc{C}}M$ the minimal non-negative
%integer $t$ such that $M^{\gg}\bigcap \mc{C}=M^{>t} \bigcap
%\mc{C}$. The notion $\mr{rAb}|_{\mc{C}}M$ is defined similarly. We
%denote $\mr{glAb}{\mc{C}}:=sup\{\mr{lAb}M\mid M\in\mc{C}\}$ and
%$\mr{glAb}_{\mc{C}}{\mc{C}}:=sup\{\mr{lAb}_{\mc{C}}M\mid
%M\in\mc{C}\}$. The notions $\mr{grAb}{\mc{C}}$ and
%$\mr{grAb}_{\mc{C}}{\mc{C}}$ are defined similarly.

%\sskip \

\bg{Rem}\label{remdef}%
$\mr{(1)}$ $\mr{LAb}M$ is just the minimal bound on the vanishing
of $\mr{Ext}_R(M,-)$ in $\mr{[\ref{Cow}]}$. In  $\mr{[\ref{HJ},
\ref{M}]}$, $\mr{lAb}M$ is also denoted by $e_R(M,-)$ in case
$M\in{\mr{mod}R}$.

$\mr{(2)}$ One can similarly define the Tor-Auslander bound of a
pair $(\mc{C},\mc{D})$, denoted by $\mr{tAb}(\mc{C},\mc{D})$, by
the bifunctor $\mr{Tor}_{i>m}^R(-,-)$. Then
$\mr{tAb}(M,{\mr{Mod}R})$ for a right $R$-module $M$ is just the
minimal bound on the vanishing of $\mr{Tor}^R(M,-)$ in
$\mr{[\ref{KK}]}$.

%One can also define the Tor $\mc{C}$-Auslander bound of $M$,
%denoted by $\mr{tAb}_{\mc{C}}M$, to be the minimal non-negative
%integer $m$ such that $\mr{KerTor}^R_{\gg}(M,-)\bigcap
%\mc{C}=\mr{KerTor}^R_{>m}(M,-)\bigcap \mc{C}$, or $\infty$ if it
%doesn't exist. Then $\mr{tAb}M (:=\mr{tAb}_{\mr{Mod}R}M)$ is just
%the minimal bound on the vanishing of $\mr{Tor}^R(M,-)$ in
%$\mr{[\ref{KK}]}$.
%

$\mr{(3)}$
 Obviously, $\mr{LAb}M\le \mr{pd}M$ and $\mr{RAb}M\le
\mr{id}M$ with the equality holds if the latter is finite. If
$M\in\mr{mod}R$, then $\mr{lAb}M=\mr{pd}M$ provided that
$\mr{pd}M<\infty$.

\ed{Rem}

\sskip\

%Let $\mc{C}\subseteq \mr{Mod}R$.
% The following lemma gives some important properties of
%classes $M^{\gg}$ and ${^{\gg}M}$.
 We say that $\mc{C}$ has
the two-out-of-three property provided that any two terms in a
short exact sequence are in $\mc{C}$ implies the third term is
also in $\mc{C}$. %
%
%
%\vskip 15pt
 The proof of the following lemma is easy.

%\vskip 15pt

%
\bg{Lem}\label{gp}%general property
Let $M\in \mr{Mod}R$.

$\mr{(1)}$ All subcategories $\mr{mod}R$, $M^{\gg}$,  and
${^{\gg}M}$ have the two-out-of-three property.

$\mr{(2)}$ All subcategories $\mr{mod}R$,  $M^{\gg}$,and
${^{\gg}M}$ are closed under direct summands and finite direct
sums. Moreover, $M^{\gg}$ $({^{\gg}M}, resp.)$ is closed under
arbitrary direct products $(direct\ sums, resp.)$.

$\mr{(3)}$ $M^{\gg}=(\Omega^iM)^{\gg}$, where $\Omega^i M$ denotes
an $i$-th syzygy of $M$.

$\mr{(4)}$ ${^{\gg}M}={^{\gg}(\Omega^{-i}M)}$, where $\Omega^{-i}
M$ denotes an $i$-th cosyzygy of $M$.

$\mr{(5)}$ If $R$ is an artin algebra and $M\in\mr{mod}R$, then
$\mathbf{D}(M^{\gg})= {^{\gg}(\mathbf{D}(M))}$.
\ed{Lem}%
%
%\Pf. Easily.
%%
%\ \hfill $\Box$

%***********************************
\vskip 15pt

\bg{Lem}\label{sum}%$\mc{C}\subseteq \mr{Mod}R$ and
Let $M, N\in \mr{Mod}R$ and $\mc{C}$ be a class of $R$-modules.

$(1)$ $\mr{Ab}(M\oplus N,{\mc{C}})\le max\{\mr{Ab}(M,{\mc{C}}),
\mr{Ab}(N,{\mc{C}})\}$.

%$(1')$ $\mr{lAb}_f(M\oplus N)\le max\{\mr{lAb}_fM, \mr{lAb}_fN\}$.
%%Moreover, if $M^{\gg}\subseteq N^{\gg}$ (e.g., $N$ has finite
%%projective dimension), then $\mr{lAb}(M\oplus N)= \mr{lAb}M$.

$\mr{(2)}$ $\mr{rAb}({\mc{C}},M\oplus N)\le
max\{\mr{Ab}(\mc{C},M), \mr{Ab}({\mc{C}},N)\}$.%, and
%$\mr{rAb}_f(M\oplus N)\le max\{\mr{rAb}_fM, \mr{rAb}_fN\}$.
%If  $^{\gg}M\subseteq {^{\gg}N}$ (e.g., $N$ has finite
%injective dimension), then $\mr{rAb}(M\oplus N)= \mr{rAb}M$.

$\mr{(3)}$ If $R$ is an artin algebra and $M\in\mr{mod}R$, then
$\mr{lAb}M=\mr{rAb}(\mathbf{D}(M))$.
\ed{Lem}%

\Pf. (1) Clearly we can assume that $k:= max\{\mr{Ab}(M,{\mc{C}}),
\mr{Ab}(N,{\mc{C}})\}<\infty$. Now note that $(M\oplus
N)^{\gg}\bigcap\mc{C}=M^{\gg}\bigcap
N^{\gg}\bigcap\mc{C}=M^{>k}\bigcap N^{>k}\bigcap\mc{C}=(M\oplus
N)^{>k}\bigcap\mc{C}$, so we have that $\mr{Ab}(M\oplus
N,{\mc{C}})\le k$. %The finitely generated version is proved similarly.
%
%If $M^{\gg}\subseteq N^{\gg}$, then $(M\oplus N)^{\gg
%0}=M^{\gg}\bigcap N^{\gg}=M^{>\mr{lAb}M}=(M\oplus
%N)^{>\mr{lAb}M}$, i.e., $\mr{lAb}(M\oplus N)\le \mr{lAb}M$.

The proof of (2) is dual to that of (1). The proof of (3) follows
from Lemma 2.3(5).% and the involved definitions.
%
%%(2) Dually to (1).
%%
%%(3) By Lemma \ref{gp} (5)  and the involved definitions.
%
\  \hfill $\Box$

%
% **************************************
%
%Let $M,N\subseteq \mr{mod}R$. The pair $(M,N)$ is called perfect
%if $M^{\gg}=N$ and ${^{\gg}N}=M$. The perfect pair is then
%said to be generated (resp., cogenerated) by $\mc{S}$ if $M={^{\gg
%0}\mc{S}}$ (resp., $N=\mc{S}^{\gg}$).
%
%**************************************
% \bg{Lem}\label{perfect}%
%Every class $M$ cogenerates a perfect pair $({^{\gg}(M^{\gg
%0})},M^{\gg})$ and generates a perfect pair $({^{\gg
%0}M},({^{\gg}M})^{\gg})$. %%
%\ed{Lem}%
%
%\Pf.
%
%**************************************
% \bg{Cor}\label{bjcor}%
% If $(M,N)$ is perfect, then $\mr{lAb}M=
%\mr{rAb}N$. %%
%\ed{Cor}%
%
%\Pf.
%
% **************************************
%
%For example, both pairs $(\mr{mod}R,\mr(I)^{<\infty})$ and
%$(\mr(P)^{<\infty},\mr{mod}R)$ are perfect, thus we have that
%$\mr{lAb}(\mr{mod}R)= \mr{rAb}(\mr(I)^{<\infty}) (=fin.injdim)$
%and $\mr{lAb}(\mr(P)^{<\infty})= \mr{rAb}(\mr{mod}R)
%(=fin.projdim)$.
%

%
%%
%%**************************************
\zskip\

%
%**************************************
%
\bg{Pro}\label{gdbj} Let ${\mc{C}}, {\mc{D}}, {\mc{E}}, {\mc{F}}$
be four classes of $R$-modules. If ${\mc{C}}\subseteq {\mc{D}}$
and ${\mc{E}}\subseteq {\mc{F}}$, then
$\mr{Ab}({\mc{C}},\mc{E})\le \mr{Ab}(\mc{C},{\mc{F}})\le
\mr{Ab}(\mc{D},{\mc{F}})$.
\ed{Pro}

\Pf. We may assume that $0\le\mr{Ab}(\mc{C},{\mc{F}})=t<\infty$.
Then $\mr{Ext}_R^{>t}(C,F)=0$ for any $C\in\mc{C}$ and
$F\in\mc{F}$ with $\mr{Ext}_R^{\gg}(C,F)=0$, by the definition. In
particular, since ${\mc{E}}\subseteq {\mc{F}}$, we have that
$\mr{Ext}_R^{>t}(C,E)=0$ for any $C\in\mc{C}$ and $E\in\mc{E}$
with $\mr{Ext}_R^{\gg}(C,E)=0$. This shows that
$\mr{Ab}({\mc{C}},\mc{E})\le t$. The remaining part is proved
similarly.
\ \hfill $\Box$

%
%**************************************
\zskip\

Let $\mc{C},\mc{D}$ be two classes of $R$-modules. We have the
following result.
%
%\vskip 15pt

\bg{Lem}\label{GD}%
$\mr{Ab}(\mc{C},\mc{D})=sup\ \{\mr{Ab}(M,{\mc{D}})\mid
M\in\mc{C}\}$\\
\indent \hskip 95pt $=sup\ \{\mr{Ab}({\mc{C}},N)\mid N\in\mc{D}\}$\\
\indent \hskip 95pt $=sup\ \{\mr{Ab}(C,D)\mid C\in\mc{C},
D\in\mc{D}\}$.
\ed{Lem}%

\Pf. By Proposition \ref{gdbj}, we have that

$sup\ \{\mr{Ab}(C,D)\mid C\in\mc{C}, D\in\mc{D}\}\le sup\
\{\mr{Ab}(M,{\mc{D}})\mid M\in\mc{C}\}\le \mr{Ab}(\mc{C},\mc{D}).$

\noindent To prove the other part, we may assume that $0\le sup\
\{\mr{Ab}(C,D)\mid C\in\mc{C}, D\in\mc{D}\}=t<\infty$. Take any
$C\in\mc{C}$ and $D\in\mc{D}$ with $\mr{Ext}^{\gg}_R(C,D)=0$, then
we easily see that $\mr{Ext}^{>t}_R(C,D)=0$. Hence, we have that
$\mr{Ab}(\mc{C},\mc{D})\le t$.
\  \ \hfill$\Box$

%
%%Let $\mc{C}$ be a subclass of $\mr{Mod}R$. now introduce the following notions.
%Let $\mc{C},\mc{D}$ be two classes of $R$-modules. %\subseteq \mr{Mod}R$.
%We set $\mr{glAb}_{\mc{C}}\mc{D}:=sup\ \{\mr{lAb}_{\mc{C}}M\mid
%M\in\mc{D}\}$ and similarly $\mr{grAb}_{\mc{C}}{\mc{D}}:=sup\
%\{\mr{rAb}_{\mc{C}}M\mid M\in\mc{D}\}$. If $\mc{C}=\mr{Mod}R$,
%then we drop $\mc{C}$ in the notions.
%
%
%%The notions $\mr{grAb}{\mc{C}}$ and
%%$\mr{grAb}_{\mc{C}}{\mc{C}}$ are defined similarly.
%%
%
%We have the following result.% gives an important relation between
%%$\mr{glAb}_{\mc{C}}{\mc{C}}$ and $\mr{grAb}_{\mc{C}}{\mc{C}}$.
%
%%
\vskip 15pt
%
%\bg{Lem}\label{GD}%
%$\mr{glAb}_{\mc{C}}{\mc{C}}= \mr{grAb}_{\mc{C}}{\mc{C}}$.
%%
%\ed{Lem}%
%
%\Pf. Assume that $\mr{grAb}_{\mc{C}}{\mc{C}}=t<\infty$. Then for
%any $M\in {\mc{C}}$ and $N\in M^{\gg}\bigcap {\mc{C}}$, we have
%that $M\in {^{\gg}N\bigcap {\mc{C}}}= {^{>t}N\bigcap{\mc{C}}}$ by
%the definition of $\mr{grAb}_{\mc{C}}{\mc{C}}$. Hence $N\in
%M^{>t}\bigcap {\mc{C}}$. It follows that
%$M^{\gg}\bigcap{\mc{C}}=M^{>t}\bigcap{\mc{C}}$. Thus,
%$\mr{lAb}_{\mc{C}}M\le t$ and $\mr{glAb}_{\mc{C}}{\mc{C}}\le t$. % by
%%the involved definitions.
%%
%Similarly, we have also $\mr{grAb}_{\mc{C}}{\mc{C}}\le
%\mr{glAb}_{\mc{C}}{\mc{C}}$.
%%
%\  \ \hfill$\Box$
%

%
%**************************************
%

%\vskip 15pt

We call $\mr{Ab}({\mc{C}},{\mc{C}})$ the {\it global Auslander
bound} of the class ${\mc{C}}$. We denote by $\mr{GAb}R$ the
global Auslander bound of $\mr{Mod}R$ and by $\mr{gAb}R$ the
global Auslander bound of $\mr{mod}R$. Note that $\mr{gAb}R$ is
just the Ext-index of $R$ in   [\ref{HJ}, \ref{M}].
%\sskip\

If $R$ is an artin algebra, then there is a duality $\mathbf{D}$
between $\mr{mod}R$ and $\mr{mod}R^o$. Hence we can easily obtain
that $\mr{gAb}R=\mr{gAb}R^o$ in this case.
%
%
%\vskip 15pt
%

%\vskip 15pt
%

%
%**************************************
%%
%\bg{Pro}\label{lim}%
%Assume that $R$ is a left coherent ring and $M\in\mr{mod}R$, Then
%$\mr{lAb}M=\mr{lAb}_{\mr{mod}R}M$. Consequently,
%$\mr{gAb}R=\mr{glAb}(\mr{mod}R)$.
%%
%\ed{Pro}
%
%\Pf. Assume that $m=\mr{lAb}_{\mr{mod}R}M<\infty$. Since $R$ is
%left coherent, every $R$-module $N$ is a direct limit of modules
%$N_i$ in $\mr{mod}R$, i.e., $N=\mr{lim}_{\to}N_i$. Now assume that
%$N\in M^{\gg}$, then we have that, for some integer $t_N$ and
%all $j>t_N$,
%$0=\mr{Ext}_R^j(M,N)=\mr{Ext}_R^j(M,\mr{lim}_{\to}N_i)$. Note
%$M\in\mr{mod}R$, so $0=\mr{Ext}_R^j(M,\mr{lim}_{\to}N_i)\simeq
%\mr{lim}_{\to}\mr{Ext}_R^j(M,N_i)$.
%
%%%%%%%%%%%%%%%下面结论有问题（极限为零不能说每个为零）%%%%%%%%%%%%%%%%%%%%%%%%%%%%%
%It follows that $\mr{Ext}_R^j(M,N_i)=0$ for all $j>t_N$.
%%%%%%%%%%%%%%%%%%%%%%%%%%%%%%%%%%%%%%%%%%%%
%
%
%Hence $N_i\in M^{\gg}\bigcap\mr{mod}R= M^{>m}\bigcap\mr{mod}R$.
%Then we obtain that
%$\mr{Ext}_R^j(M,N)=\mr{lim}_{\to}\mr{Ext}_R^j(M,N_i)=0$ for all
%$j>m$. Therefore, $\mr{lAb}M\le m$.
%
%\ \hfill $\Box$

%
%**************************************
%
%\bskip\

To calculate the Auslander bound of a module, it is enough to
calculate the Auslander bound of its syzygies, as the follow
result shows.

\bg{Lem}\label{syg}%Let $\mc{C}\subseteq \mr{Mod}R$.
Let $M\in \mr{Mod}R$ and $\mc{C}$ be a class of $R$-modules.

$\mr{(1)}$ $0\le\mr{Ab}(M,{\mc{C}})\le m$ if and only if
$\mr{Ab}(\Omega^mM,{\mc{C}})=0$, where $\Omega^mM$ denotes an
$m$-th syzygy of $M$.

$\mr{(2)}$ $0\le\mr{Ab}({\mc{C}},M)\le m$ if and only if
$\mr{Ab}({\mc{C}},\Omega^{-m}M)=0$, where $\Omega^{-m}M$ denotes
an $m$-th cosyzygy of $M$.
\ed{Lem}%

\Pf. %We still prove (1) and leave (2) to the readers.
(1) Since $M^{\gg}=(\Omega^mM)^{\gg}$ by Lemma \ref{gp} (3) and
$M^{>m}=(\Omega^mM)^{\ge 1}$ by dimension shifting, we see that
$0\le \mr{Ab}(M,{\mc{C}})\le m \Longleftrightarrow \O\neq
M^{\gg}\bigcap\mc{C}\subseteq M^{>m}\bigcap\mc{C}
$\\
\indent\hskip 175pt$\Longleftrightarrow
\O\neq(\Omega^mM)^{\gg}\bigcap\mc{C}\subseteq (\Omega^mM)^{\ge
1}\bigcap\mc{C}$\\
\indent\hskip 175pt$\Longleftrightarrow
\mr{Ab}(\Omega^mM,{\mc{C}})=0$.

The proof of (2) is dual to that of (1).
\ \hfill $\Box$

%%%%%%%%%%%%%%%%%%%%%%%%%%%%%%%%%%%%%%%%%%%%
\vskip 15pt
By the above lemma, we easily obtain the following result.% have % and the involved definitions

%
%**************************************
%
\bg{Pro}\label{gdkh}%
The following are equivalent for a class $\mc{C}$ such that
$\mr{Ab}(\mc{C},\mc{C})\neq -1$.

$\mr{(1)}$ The global Auslander bound of $\mc{C}$ is not more than
$n$.

$\mr{(2)}$ $\mr{Ab}(\Omega^n\mc{C},{\mc{C}})=0$, where
$\Omega^n\mc{C}$ denotes the class of all $n$-th syzygies of
$R$-modules in $\mc{C}$.

$\mr{(3)}$ $\mr{Ab}({\mc{C}},\Omega^{-n}\mc{C})=0$, where
$\Omega^{-n}\mc{C}$ denotes the class of all $n$-th cosyzygies of
$R$-modules in $\mc{C}$.

\ed{Pro}
%
%**************************************

%
%\Pf. Since $R$ is an artin algebra, there is a duality
%$\mathbf{D}$ between $\mr{mod}R$ and $\mr{mod}R^o$. Thus, we
%obtain that  $\mr{gAb}R\le n$ if and only if
%$\mr{Ext}_R^{i>n}(M,N)=0$ for all $M,N\in\mr{mod}R$ such that
%$\mr{Ext}_R^i(M,N)=0$ for all sufficiently large $i$, if and only
%if $\mr{Ext}_R^{i>n}(\mathbf{D}M,\mathbf{D}N)=0$ if and only if
%$\mr{gAb}R^o\le n$. \
%%
%\ \hfill $\Box$
%
%
%**************************************
%
%

%\vskip 15pt

In some cases, the Auslander bound of a module can be tested by
special modules, as the following theorem shows.

\bg{Th}\label{mth}%
Let $M\in \mr{Mod}R$.

$\mr{(1)}$ If $\mr{lAb}M<\infty$ and $R\in M^{\gg}$, then
$\mr{lAb}M=min\{t| R\in M^{>t}\}$.

$\mr{(2)}$ If $\mr{LAb}M<\infty$ and $R^{(\kappa)}\in M^{\gg}$ for
all cardinals $\kappa$, then $\mr{LAb}M=min\{t| R^{(\kappa)}\in
M^{>t}$ for all cardinals $\kappa\}$.

$\mr{(3)}$ If $\mr{RAb}M<\infty$ and $I\in {^{\gg}M}$ for all
injective $R$-modules $I$, then $\mr{RAb}M=min\{t| I\in {^{>t}M}$
for all injective $R$-modules $I\}$.
\ed{Th}%

\Pf. (1) Assume that $\mr{lAb}M=m<\infty$. Let $t=min\{t| R\in
M^{>t}\}$. Then $t<\infty$, since $R\in M^{\gg}$ by the
assumption.

Note that $R\in\mr{mod}R$, so we have that $R\in M^{>t}\bigcap
{\mr{mod}R}\subseteq M^{\gg}\bigcap {\mr{mod}R}$. Now take any
$N\in M^{\gg}\bigcap {\mr{mod}R}$ and any projective resolution of
$N$: $\cdots\to P_1\to P_0\to N\to 0$ with each $P_i$ finitely
generated projective. Then we have all $P_i\in M^{>t}\bigcap
{\mr{mod}R}$, and hence all $\Omega^iN\in M^{\gg}\bigcap
{\mr{mod}R}$, by Lemma \ref{gp} (1). Therefore, for all $i>t$, we
obtain that $\mr{Ext}_R^i(M,N)\simeq \mr{Ext}_R^{i+1}(M,\Omega
N)\simeq\cdots\simeq \mr{Ext}_R^{i+m}(M,\Omega^mN)=0$, by
dimension shifting and the definition of $\mr{lAb}M$. It follows
that $N\in M^{>t}\bigcap {\mr{mod}R}$. Consequently,
$M^{>t}\bigcap {\mr{mod}R}=M^{\gg}\bigcap {\mr{mod}R}$, that is,
$\mr{lAb}M\le t$.

On the other hand, since $R\in M^{\gg}\bigcap {\mr{mod}R}$ by
assumptions, we have that $R\in M^{>m}$ by the definition of
$\mr{lAb}M$. It follows that $t\le m$ too. Hence the conclusion
follows.

%\vskip 10pt

The proof of (2) is similar as (1) and the proof of (3) is dual to
(2).
\  \ \hfill $\Box$

%**************************************
%
%**************************************
%

\vskip 15pt

Immediately, we obtain the following corollary [\ref{M}, Corollary
3.3].

\bg{Cor}\label{Cor-id}%
Assume that $\mr{id}R<\infty$.

$(1)$ If $\mr{lAb}M<\infty$ for every
$M\in\mr{mod}R$, then $\mr{gAb}R=\mr{id}R$. %In the case,
%$\mr{pd}M<\infty$ if and only if $M\in M^{\gg}$, for any
%$M\in\mr{mod}R$.

$(2)$ If $R$ is left noetherian and $\mr{LAb}M<\infty$ for every
$M\in\mr{Mod}R$, then $\mr{GAb}R=\mr{id}R$. %In the case,
%$\mr{pd}M<\infty$ if and only if $M\in M^{\gg}$, for any
%$M\in\mr{Mod}R$.
%
\ed{Cor}%

%
%**************************************
%
%\sskip\

We note that assumptions in Theorem \ref{mth} (1) can not be
removed. For example, let $R$ be an artin algebra of finite
representation type with $\mr{id}R=\infty$. Then it is easy to see
that $\mr{gAb}R<\infty$. However, it is obvious that there are
modules $M\in\mr{mod}R$ such that $min\{t| R\in M^{>t}\}=\infty$.
Thus the condition $R\in M^{\gg}$ is needed. Now let $R$ be a
Gorenstein ring with $\mr{gAb}R=\infty$ (such rings exist by
[\ref{JS}]). Then there are modules $M\in\mr{mod}R$ such that
$\mr{lAb}M=\infty$ by Corollary \ref{Cor-id} (1). However, it is
easy to see that $min\{t| R\in M^{>t}\}\le\mr{id}R<\infty$ for any
$M\in\mr{mod}R$. So the condition that $\mr{lAb}M<\infty$ is also
needed.

We remark that it is an openh question whether
$\mr{GAb}R=\mr{gAb}R$ if $R$ is left noetherian.

The following theorem is our main result which relates Auslander
bounds to Auslander-Reiten conjecture as claimed in the
introduction.
%
%**************************************
%
\vskip 15pt

\bg{Th}\label{proself}% $R$ be a ring and
Let $M\in \mr{Mod}R$.

$\mr{(1)}$ Assume that $M\oplus R\in M^{\gg}\bigcap\mr{mod}R$,
then $\mr{lAb}M=\mr{pd}M$.

$\mr{(2)}$ Assume that $M\oplus R^{(\kappa)}\in M^{\gg}$ for any
cardinal $\kappa$, then $\mr{LAb}M=\mr{pd}M$.

$\mr{(3)}$ Assume that $M,I\in {^{\gg}M}$ for any injective
$R$-module $I$, then $\mr{RAb}M=\mr{id}M$.
\ed{Th}%

\Pf. (1) Clearly we need only prove that $\mr{pd}M\le \mr{lAb}M$.

We can assume that $\mr{lAb}M=m<\infty$. Since $M\in\mr{mod}R$, we
can take a projective resolution of $M$: $\cdots\to P_1\to P_0\to
M\to 0$ with each $P_i$ finitely generated projective. Then each
$\Omega^iM\in M^{\gg}\bigcap\mr{mod}R$ by Lemma \ref{gp} (1), as
$M\oplus R\in M^{\gg}\bigcap\mr{mod}R$. It follows that
$\Omega^iM\in M^{>m}$ for each $i$, by the definition of
$\mr{lAb}M$. Now by applying the functor
$\mr{Hom}_R(-,\Omega^{m+1}M)$ to the exact sequence $0\to
\Omega^mM\to P_{m-1}\to\cdots\to P_0\to M\to 0$, we obtain that
$\mr{Ext}_R^1(\Omega^mM,\Omega^{m+1}M)\simeq
\mr{Ext}_R^2(\Omega^{m-1}M,\Omega^{m+1}M)\simeq\cdots\simeq
\mr{Ext}_R^{m+1}(M,\Omega^{m+1}M)$ by dimension shifting. The
latter is 0 since $\Omega^{m+1}M\in M^{>m}$ by the above argument.
It follows that the exact sequence $0\to \Omega^{m+1}M\to P_m\to
\Omega^mM\to 0$ splits, and consequently, $\mr{pd}M\le m$. \

The proof of (2) is similar as (1) and the proof of (3) is dual to
(2).
%
%(2) Similarly as (1).
%
%(3) Dually to (2).
%
\ \hfill$\Box$

%
%**************************************
%

%
%
\vskip 30pt
% ----------------------------------------------------------------------
% ----------------------------------------------------------------------

\section{Homological conjectures}
%{\noindent \Large\bf Introduction}
%
%\vskip 15pt
%
%
\hskip 15pt

As pointed out in [\ref{JS}, \ref{S}], Auslander's conjecture
fails for artin algebras in general. However, we can consider a
finitistic version of Auslander's conjecture. Let $R$ be a ring.
We set

\hskip 15pt $\bullet\ \ \mr{FLAb}(R):=sup\ \{\mr{LAb}M\mid
\mr{LAb}M<\infty\}$, and

\hskip 15pt $\bullet\ \ \mr{flAb}(R):=sup\ \{\mr{lAb}M\mid
M\in\mr{mod}R$ and $\mr{lAb}M<\infty\}$.

Similarly, we have notions $\mr{FRAb}(R)$ and $\mr{frAb}(R)$
defined by right Auslander bounds.  Note that
$\mr{flAb}(R)=\mr{frAb}(R^o)$ in case that $R$ is an artin
algebra.

Now we formulate the following conjecture.

\bg{verse}
%\noindent

 {\bf Finitistic Auslander Conjecture} (FAC for short):
\newline
\hskip 15pt $\bullet\ $(lFAC):\ \ \ {\it $\mr{flAb}(R)<\infty$ for
every artin algebra
$R$,} or dually %(and equivalently by the duality in artin
%algebras),
\newline
\hskip 15pt $\bullet\ $(rFAC):\ \ \ {\it $\mr{frAb}(R)<\infty$ for
every artin algebra $R$.}

\ed{verse}

It is easy to see that the finitistic Auslander conjecture implies
the finitistic dimension conjecture for artin algebras by Remark
\ref{remdef} (3).

It is also clear that $\mr{flAb}R=\mr{frAb}R<\infty$ if
$\mr{gAb}R<\infty$, by Lemma \ref{GD}. Similarly,
$\mr{FLAb}R=\mr{FRAb}R<\infty$ if $\mr{GAb}R<\infty$. For example,
every group algebra $kG$ with $k$ a field and $G$ finite has the
property $\mr{GAb}(kG)<\infty$, see [\ref{BCR}, Theorem 2.4] and
[\ref{CH}, Appendix A].
%
%It is also interesting to know whether $\mr{fLAb}(R)=\mr{FLAb}(R)$
%for every artin algebra $R$.

%%
%\vskip 15pt
%%
%
%
%
%\bg{Pro}\label{fac-fdc}%
%Let $R$ be a ring.
%
%$\mr{(1)}$ $\mr{fLAb}(R)<\infty$ implies $\mr{fPD}(R)<\infty$.
%
%$\mr{(1')}$ $\mr{FLAb}(R)<\infty$ implies $\mr{FPD}(R)<\infty$.
%
%$\mr{(2)}$ $\mr{fRAb}(R)<\infty$ implies $\mr{fID}(R)<\infty$.
%
%$\mr{(2')}$ $\mr{FRAb}(R)<\infty$ implies $\mr{FID}(R)<\infty$ .
%%
%%$\mr{(3)}$ If $R$ is left noetherian with $\mr{gAb}R<\infty$, then
%%$\mr{fPD}(R)<\infty$ and $\mr{fID}(R)<\infty$.
%%
%%$\mr{(3')}$ If $R$ is left noetherian with $\mr{gAb}R<\infty$,
%%then $\mr{fPD}(R)<\infty$ and $\mr{fID}(R)<\infty$.
%%
%\ed{Pro}%
%
%\Pf. Directly by Remark \ref{remdef} (3).
%%
%\ \ \hfill$\Box$

%\vskip 15pt

As  the finitistic dimension conjecture fails for commutative
noetherian rings in general, the conjecture (FAC) fails in the
case, too. Moreover, it is pointed out in [\ref{CH}] that there is
a commutative notherian ring $R$ with infinite Krull dimension
such that $\mr{lAb}M<\infty$ for every $M\in\mr{mod}R$ but
$\mr{gAb}R=\infty$.

%However, it is interesting to check (FAC) in case of artin
%algebras. In particular,

It is unknown whether $\mr{gAb}R<\infty$ if $R$ is an artin
algebra such that $\mr{lAb}M<\infty$ for every $M\in\mr{mod}R$.

The following result gives a partial answer to the conjecture
(FAC).

\vskip 15pt

\bg{Pro}\label{id-fac}%
Let $R$ be a ring.

$(1)$ If $\mr{id}R<\infty$, then $\mr{flAb}(R)<\infty$.

$(2)$ If $R$ is left noetherian and $\mr{id}R<\infty$, then
$\mr{FLAb}(R)<\infty$.
\ed{Pro}%

\Pf. (1) Indeed, we have that $\mr{lAb}M=min\{t| R\in M^{>t}\}\le
\mr{id}R$ provided $\mr{lAb}M<\infty$ and $M\in\mr{mod}R$, by
Theorem \ref{mth}.

(2) If  $R$ is left noetherian and $\mr{id}R<\infty$, then
$\mr{id}R^{(\kappa)}<\infty$ for any cardinal $\kappa$. Now the
conclusion follows from  Theorem \ref{mth} again.
\ \hfill$\Box$

%
%**************************************
%

%
\vskip 15pt

Now we turn to other related homological conjectures.

Let $R$ be a ring and $T\in \mr{Mod}R$ with $S=\mr{End}_RT$.
Recall from [\ref{W}] that  $T$ is Wakamatsu-tilting if it
satisfies

\begin{verse}
(1) $T\in \mr{mod}R$ and $T_S\in\mr{mod}S^o$,\\
(2) $R\simeq \mr{End}(T_S)$, and\\
(3) $\mr{Ext}_R^i(T,T)=0=\mr{Ext}_{S^o}^i(T,T)$ for all $i>0$.
\end{verse}

Equivalently, as shown in [\ref{W}],  $T$ is Wakamatsu-tilting if

\begin{verse}
(W1) $T\in \mr{mod}R$,\\
(W2) $\mr{Ext}_R^i(T,T)=0$  for all $i>0$, and\\
(W3) there is an exact sequence $0\to R\to^{f_0} T_0\to^{f_1}
T_1\to^{f_2}\cdots$ with each $T_i\in\mr{add}T$ and each
$\mr{Im}f_i\in {^{>0}T}$, for all $i\ge 0$.
\end{verse}

It is clear that $T$ is Wakamatsu-tilting if and only if $T_S$ is
Wakamatsu-tilting.

Recall also that $T$ is tilting [\ref{AR}, \ref{My}] if it
satisfies

\begin{verse}
(T1) $T\in \mr{mod}R$ and $\mr{pd}T<\infty$,\\
(T2) $\mr{Ext}_R^i(T,T)=0$ for all $i>0$, and\\
(T3) there is an exact sequence $0\to R\to T_0\to\cdots \to T_n\to
0$ with each $T_i\in\mr{add}T$, for some integer $n$.
\end{verse}

 We note that
$T$ is tilting if and only if $T_S$ is tilting, where
$S=\mr{End}_RT$ [\ref{My}]. The following conjecture is cited from
[\ref{MR}].

% If $R$ is
%an artin algebra, then $T$ is $n$-cotilting if it satisfies (C1)
%$\mr{id}T\le n$, (C2) $\mr{Ext}_R^i(T,T)=0$ for all $i>0$, and
%(C3) there is an exact sequence $0\to T_n\to\cdots\to T_0\to
%\mathbf{D}(R_R)\to 0$ with each $T_i\in\mr{add}T$. Obviously, $T$
%is tilting if and only if $(\mathbf{D}T)_R$ is cotilting in case
%$R$ is an artin algebra.

%
\bg{verse}
%
%\noindent

 {\bf Wakamatsu Tilting Conjecture}
%
%\newline\hskip 15pt
(WTC for short):  {\it Every Wakamatsu-tilting module of finite
projective dimension is tilting.} %Or dually (and equivalently by
%the duality in artin algebras),
%
%%
%(IWTC)  {\it Every Wakamatsu-tilting module of finite injective
%dimension is tilting.}

%
\ed{verse}
%
%
%
%It is clear that (PWTC) holds for an artin algebra $R$ if and only
%if (IWTC) holds for $R^o$.

It is pointed out in [\ref{MR}] that, if the finitistic dimension
conjecture holds for a ring $R$, then the conjecture (WTC) holds
for $R$. We have an equivalent version of the conjecture (WTC)
(and so we denote this conjecture by EWTC, where E means
equivalent).

\vskip 15pt

\bg{Pro}\label{EWTC}%
The conjecture $\mr{(WTC)}$ holds for all rings if and only if the
following conjecture $\mr{(EWTC)}$ holds for all rings $R$.

\bg{verse}

$\mr{(EWTC)}:$ An $R$-module $T\in\mr{mod}R$ is tilting if it
satisfies conditions $\mr{(T2)}$ and $\mr{(T3)}$ in the definition
of tilting modules.\ed{verse}
\ed{Pro}%

\Pf. Let $R$ be a ring and $T\in \mr{mod}R$ with $S=\mr{End}_RT$.

$\mr{(EWTC)} \Rightarrow \mr{(WTC)}$: Assume that $T$ is a
Wakamatsu-tilting $R$-module with $\mr{pd}T<\infty$. Then we have
an exact sequence $0\to P_n\to\cdots\to P_0\to T\to 0$ for some
$n$. Applying the functor $\mr{Hom}_R(-,T)$, we obtain an induced
exact sequence $0\to S\to T_0\to\cdots\to T_n\to 0$ in
$\mr{mod}S^o$, since $\mr{Ext}_R^i(T,T)=0$ for all $i>0$. Since
$T$ is a Wakamatsu-tilting $R$-module,  $T$ is also a
Wakamatsu-tilting $S^o$-module and so $\mr{Ext}_{S^o}^i(T,T)=0$
for all $i>0$. Hence we get that $T_S$ is tilting provided that
(EWTC) holds for $S^o$. Consequently, $T$ is also a tilting
$R$-module.

$\mr{(WTC)} \Rightarrow \mr{(EWTC)}$: If (T2) and (T3) in the
definition of tilting modules holds for $T\in\mr{mod}R$, then $T$
is Wakamatsu-tilting. Moreover, by applying the functor
$\mr{Hom}_R(-,T)$ to the exact sequence in (T3), we easily see
that $T_S$ is Wakamatsu-tilting with finite projective dimension.
It follows that $T_S$ is tilting  provided that (WTC) holds for
$S^o$. Now by the left-right symmetry, we get that $T$ is tilting.
\  \ \hfill$\Box$

\vskip 15pt

If we specify $n=0$ in the condition (T3), then $R\in \mr{add}T$.
In this case, the conjecture (EWTC) is just the Auslander-Reiten
conjecture.
%
%%
%\vskip 15pt
%%
%
%\noindent {\bf Auslander-Reiten Conjecture:}
%%
%(ARC)  {\it Every selforthogonal generator possessing finitely
%generated projective resolutions is projective.}
%
%%
%\vskip 15pt
%%

 The following result gives a partial answer to the conjecture (EWTC), which
extends [\ref{CH}, Theorem A].

\vskip 15pt

\bg{Pro}\label{Ewtch}%
Let $R$ be a ring. If $\mr{lAb}M<\infty$ for every
$M\in\mr{mod}R$, then the conjecture $\mr{(EWTC)}$ holds for $R$.
%In particular, $\mr{(ARC)}$ holds for $R$.
%
\ed{Pro}%

\Pf. Assume that $T\in\mr{mod}R$ satisfies the conditions (T2) and
(T3). Then we easily obtain that $T\oplus R\in T^{\gg}$. Now by
the assumption and Theorem \ref{proself}, we get that
$\mr{pd}T=\mr{lAb}T<\infty$. It follows that $T$ is
tilting.% by the involved definition.
\  \ \hfill$\Box$

%%
%% **************************************
%%
%\vskip 15pt
%
%The conjecture (EWTC) has an infinitely generated version. Let $R$
%be a ring and $T\in\mr{Mod}R$. Recall from [\ref{AhC}] that $T$ is
%infinitely generated tilting if it satisfies (T1')
%$\mr{pd}T<\infty$, (T2') $\mr{Ext}_R^i(T,T^{(\kappa)})=0$ for all
%cardinals $\kappa$ and all $i>0$, and (T3') there is an exact
%sequence $0\to R\to T_0\to\cdots \to T_n\to 0$ with each
%$T_i\in\mr{Add}T$. Similarly as the above result, we also have the
%following.
%
%
%\vskip 15pt
%
%
%
%\bg{Pro}\label{EWTCI}%
%Let $R$ be a ring and $T\in\mr{Mod}R$. If $\mr{lAb}M<\infty$ for
%every $M\in\mr{Mod}R$, then $T$ is infinitely generated tilting if
%and only if it satisfies $\mr{(T2')}$ and $\mr{(T3')}$ in the
%definiton. In particular, every selforthogonal generator is
%projective.
%%
%\ed{Pro}%
%
%\Pf. Similarly as Proposition \ref{Ewtch}.
%%
%\  \ \hfill$\Box$

%
% **************************************
%
\zskip\

%Another special case of (WTC) is as follows (indeed, restricting
%to the trivial Wakamatsu-tilting module $\mathbf{D}R$).

We now consider another homological conjecture.
%
%\noindent

\bg{verse} {\bf Gorenstein Symmetry Conjecture} :
{\it Let $R$ be an artin algebra. Then $\mr{id}R<\infty$ if and
only if $\mr{id}(R_R)<\infty$.}

\ed{verse}

Gorenstein Symmetry conjecture clearly makes sense for any ring.
It was proved in [\ref{CH}] that if $R$ is a two-sided noetherian
ring such that $\mr{lAb}M<\infty$ for every $M\in\mr{mod}R$ and
$\mr{lAb}N<\infty$ for every $N\in\mr{mod}R^o$, and (1) $R$ is an
artin algebra, or (2) $R$ has a dualizing complex, then
$\mr{id}R<\infty$ if and only if $\mr{id}_{R^o}R<\infty$ (whence,
$\mr{id}R = \mr{id}_{R^o}R$ by [\ref{I}]). The following result
also gives a similar answer to the Gorenstein Symmetry conjecture.
Note that we do not know if $\mr{LAb}M<\infty$ for every
$M\in\mr{Mod}R$ provided that $\mr{lAb}M<\infty$ for every
$M\in\mr{mod}R$, even when $R$ is an artin algebra. We do not know
whether $\mr{LAb}M<\infty$ for every $M\in\mr{Mod}R$ implies that
$\mr{RAb}M<\infty$ for every $M\in\mr{Mod}R$.

%
%**************************************
%

\vskip 15pt

\bg{Pro}\label{gsc}%
Let $R$ be a two-sided noetherian ring. Assume that

\begin{verse}
$(1)$ $\mr{LAb}M<\infty$ for every $M\in\mr{Mod}R$ and every
$M\in\mr{Mod}R^o$, or\\
$(2)$ $R$ has a dualizing complex and $\mr{rAb}M<\infty$ for every
$M\in\mr{mod}R$ and every $M\in\mr{mod}R^o$.
\end{verse}

Then $\mr{id}R<\infty$ if and only if $\mr{id}R^o<\infty$.
\ed{Pro}%

\Pf. Assume (1) holds. If $\mr{id}R<\infty$, then
$\mr{GAb}R<\infty$ by assumptions and Corollary \ref{Cor-id} (2).
It follows that $\mr{FRAb}(R)<\infty$. Using the fact that
$\mr{RAb}M= \mr{id}M$ for any $R$-module $M$ with
$\mr{id}M<\infty$ in Remark \ref{remdef} (3), we further obtain
that the injective version of the finitistic dimension conjecture
holds for $R$. This shows that $\mr{id}R^o<\infty$ by [\ref{FW},
Proposition 7]. Similarly, we can prove that if
$\mr{id}R^o<\infty$ then $\mr{id}R<\infty$.

Assume now (2) holds. If $\mr{id}R<\infty$, then for any
$M\in\mr{mod}R$ and any injective $R^o$-module $N$, it holds that
$\mr{Tor}_{(\mr{id}R)+1}^R(N,M)\simeq
\mr{Hom}_{R^o}(\mr{Ext}_R^{(\mr{id}R)+1}(M,R),N)=0$. Hence
$\mr{fd}_{R^o}N<\infty$. By the definition of dualizing complex
[\ref{CH}, Section 3.4], $R^o$ has a dualizing complex if so is
$R$. In this case, we have that all $R^o$-modules of finite flat
dimension have finite projective dimension, by [\ref{Jor},
Theorem]. It follows that $N\in {^{\gg}_{\ \ R^o}R}$ for any
injective $R^o$-module $N$. Now applying Theorem \ref{proself} to
the $R^o$-module $R$, we obtain that $\mr{id}_{R^o}R=\mr{rAb}R$
and the latter is finite by assumptions. Thus, we have that
$\mr{id}R^o<\infty$. The proof of the other part is also similar.
\  \ \hfill$\Box$

%
% **************************************
%

%%%%%%%%%%%%%%%%%%%%%%%%%%%%%%%%%%%%%%%%%%%%%%%%%%%%%%%%%%%%%%%%%%%%%%%%%%%%%%%%%%
\vskip 15pt

In the remaining part, we discuss a class of rings with finite
global Auslander bound.

Let $R$ be a ring and $M$ an $R$-module. Assume that $n$ is a
nonnegative integer. Following Goodearl and Zimmermann-Huisgen
[\ref{GZ}], we say that a projective resolution of $M$ is
repetitive at degree $n$ if there exists a decomposition
$\Omega^n(M) = P\oplus A_i$ such that $P$ is projective and each
$A_i$ occurs as a direct summand of \textit{infinitely many}
$\Omega^j(M)$. The repetition index of $M$, denoted $\mr{rep}(M)$,
is the least nonnegative integer $k$ such that there is a %minimal
projective resolution of $M$  which is repetitive at degree $k$
(if such a $k$ exists), or $\infty$ (otherwise). The corresponding
global repetition index is $\mr{Grep}(R)= sup\ \{\mr{rep}(M)|M\in
\mr{Mod}R\}$.

We have the following result which relates the repetition index to
the Auslander bound.

\vskip 15pt

\bg{Lem}\label{rp}%
Let $R$ be a ring and $M\in\mr{Mod}R$. If $\mr{rep}(M)=m<\infty$,
then $\mr{LAb}M\le m$.
\ed{Lem}%

\Pf. Since $\mr{rep}(M)=m<\infty$, $M$ has a projective resolution
such that $\Omega^m(M) = P\oplus A_i$ with $P$ projective and that
each $A_i$ occurs as a direct summand of infinitely many
$\Omega^j(M)$. So we have that $\Omega^mM\in\mr{add}( P\oplus
(\oplus_{i>t}\Omega^iM))$, for any $t$. Now take any $N\in
M^{\gg}$ and assume that $N\in M^{>t}$ for some $t$. It follows
that $\mr{Ext}_R^1(\Omega^mM, N)\le
\mr{Ext}_R^1(\oplus_{i>t}\Omega^iM, N)\simeq
\prod_{i>t}\mr{Ext}_R^{i+1}(M, N)=0$ by dimension shifting and the
definition of $M^{>t}$. Note that $(\Omega^mM)^{\gg}=M^{\gg}$ by
Lemma \ref{gp}, so $\mr{LAb}(\Omega^mM)=0$. It follows that
$\mr{LAb}M\le m$ by Lemma \ref{syg}.
\ \hfill$\Box$

\vskip 15pt

Consequently, we obtain the following result.

\bg{Pro}\label{uc}%
If $R$ is a ring with $\mr{Grep}R<\infty$, then
$\mr{GAb}R<\infty$. In this case, the conjecture $\mr{(EWTC)}$
holds. In particular, the Auslander-Reiten conjecture holds in
this case.
\ed{Pro}%

\Pf. By Lemma \ref{rp} and Proposition \ref{Ewtch}.
\ \hfill$\Box$

\vskip 15pt

%
%**************************************
%

Let $R$ be an artin ring. Recall that $M\in\mr{mod}R$ has a
ultimately closed projective resolution if there is some $m$ such
that the $m$-th syzygy $\Omega^mM\simeq \oplus M_i$ with each
$M_i\in\mr{add}(\Omega^{m_i}M)$ for some $m_i<m$. In this case, we
have that $\mr{rep}(M)\le m$, see [\ref{GZ}]. One defines artin
rings such that every finitely generated module has an ultimately
closed projective resolution to be of
 projective ultimately closed type.
It was proved in [\ref{AR}] that the  Auslander-Reiten conjecture
holds for artin algebras of projective ultimately closed type.
Lemma \ref{rp} and Proposition \ref{Ewtch} together also imply
that the conjecture (EWTC) and Auslander's conjecture hold for
such artin algebras.

In [\ref{GZ}], the authors studied finite dimensional algebra
$R=\mc{O}/\pi\mc{O}$, where $\mc{O}$ is a classical order over a
discrete valuation ring $D$ with uniformizing parameter $\pi$ and
residue class field $K$. The homological properties of $\mc{O}$
are to a great extent determined by those of $R$ while the latter
algebra is substantially easier to handle. In their paper, it was
shown that, if $\mr{gd}\mc{O} = d < 1$, then $\mr{Grep}R= d-1$, in
particular, the finitistic dimension of $R$ is finite. Combining
these with results in this paper and [\ref{CH}], we also know that
in case $\mc{O}$ has finite global dimension, the algebra
$\mc{O}/\pi\mc{O}$ also satisfies the Auslander conjecture, the
conjecture (FAC), the Auslander-Reiten conjecture, Gorenstein
Symmetry conjecture and the conjecture (EWTC).

\vskip 10pt

Results in this section suggest the following conjecture which
generalizes Auslander-Reiten conjecture.

\bg{verse} {\bf Generalized Auslander-Reiten Conjecture} : Let $R$
be a ring and $M\in\mr{mod}R$. If $M\oplus R\in M^{\gg}$, then
$\mr{pd}M<\infty$.\ed{verse}

%
%\vskip 10pt
%
%\noindent {\bf Conjecture B:} Let $R$ be a ring and
%$M\in\mr{mod}R$. If $M\in M^{\gg}$, then $\mr{lAb}M<\infty$.
%

\vskip 10pt

%By Theorem \ref{proself}, Conjecture B implies (GARC). Thus,

By Lemma \ref{rp}, the conjecture holds for artin algebras of
projective ultimately closed type.

%
%%%%%%%%%%%%%%%%%%%%%%%%%%%%%%%%%%%%%%%%%%%%%%%%%%%%%%%%%%%%%%%%%%%%%%%%%%%%%%%%%
%**************************致 谢**********************************************
\vskip 20pt
\begin{center}
{\bf ACKNOWLEDGEMENTS}
\end{center}

It is a pleasure to thank the referee for his/her carefully
reading and excellent suggestions.
%%%%%%%%%%%%%%%%%%%%%%%%%%%%%%%%%%%%%%%%%%%%%%%%%%%%%%%%%%%%%%%%%%%%%%%%%%%%%%%%%
%**************************参考文献**********************************************
%
\vskip 30pt\baselineskip 18pt %\baselineskip 12pt %\baselineskip 16pt
{\small %

\vskip 80pt{\it

%\noindent Address:

\noindent Jiaqun Wei

\noindent  School of Mathematics Science, Nanjing Normal
University, Nanjing 210046, China

\noindent Email: weijiaqun@njnu.edu.cn }

\end{document}